\documentclass[12pt]{amsart}
\usepackage{amsmath,amsthm,amssymb}
\usepackage[mathcal]{euscript}
\usepackage{color}
\usepackage[all]{xy}
\usepackage{epsf}

\newcommand{\blue}[1]{{\color{blue}#1}}
\newcommand{\red}[1]{{\color{red}#1}}

\newtheorem{theorem}{Theorem}[section]
\newtheorem{proposition}[theorem]{Proposition}
\newtheorem{lemma}[theorem]{Lemma}
\newtheorem{corollary}[theorem]{Corollary}

\theoremstyle{definition}

\newtheorem{remark}[theorem]{Remark}

\newcommand{\abs}[1]{\lvert#1\rvert}   

\newcommand{\st}{\mathrm{st}}
 \newcommand{\id}{\mathit{id}} 
  \newcommand{\Id}{\mathit{Id}}          
\newcommand{\Inv}{\mathrm{Inv}}
 \newcommand{\Des}{\mathrm{Des}}
 
 \newcommand{\GDes}{\mathrm{GDes}}      
 \newcommand{\Sh}{\mathrm{Sh}}
 
\newcommand{\bracket}[1]{\{#1\}}   
\newcommand{\ebracket}{\bracket{\ ,\ } }  
\newcommand{\mat}[1]{\begin{bmatrix}#1\end{bmatrix}} 

\newcommand{\map}[1]{\xrightarrow{#1}}
\newcommand{\leftmap}[1]{\xleftarrow{#1}}

\newcommand{\under}{\backslash}

\newcommand{\calQ}{\mathcal{Q}}    
  
\newcommand{\calK}{\mathcal{K}}    
\newcommand{\calY}{\mathcal{Y}}    
\newcommand{\calL}{\mathcal{L}} 
\newcommand{\calC}{\mathcal{C}} 
\newcommand{\calP}{\mathcal{P}}
\newcommand{\calS}{\mathcal{S}} 
\newcommand{\calF}{\mathcal{F}} 

\newcommand{\sfL}{\mathsf{L}} 
\newcommand{\sfD}{\mathsf{D}} 

\newcommand{\Ass}{\mathcal{A}s}
\newcommand{\Lie}{\mathcal{L}ie}
\newcommand{\Anon}{\mathcal{A}}
\newcommand{\Lnon}{\mathcal{L}}

\newcommand{\field}{\Bbbk} 

\newcommand{\setS}{\mathsf{S}}
\newcommand{\setT}{\mathsf{T}}
\newcommand{\setR}{\mathsf{R}}

\begin{document}

\title[Associative operad and weak order]{The associative operad and the weak order on the symmetric groups}
\author[M. Aguiar]{Marcelo Aguiar}
\address{Department of Mathematics\\
         Texas A\&M University\\
         College Station,  TX  77843\\
         USA}
\email{maguiar@math.tamu.edu}
\urladdr{http://www.math.tamu.edu/$\sim${}maguiar}

\author[M. Livernet]{Muriel Livernet}
\address{Institut Galil\'ee\\
Universit\'e Paris Nord\\
  93430 Villetaneuse       \\
         France}
\email{livernet@math.univ-paris13.fr}
\urladdr{http://www.math.univ-paris13.fr/$\sim${}livernet/}

\thanks{Aguiar supported in part by NSF grant DMS-0302423. 
We thank the  Institut Galil\'ee of the Universit\'e Paris 13 for a one-month invitation that led to this work.}

\keywords{Operad, coalgebra, permutation, weak Bruhat order, binary tree, composition}
\subjclass[2000]{Primary 18D50, 06A11; Secondary 06A07, 16W30}

\date{October 17, 2006}

\begin{abstract} The associative operad  is a certain algebraic structure on the sequence of group algebras of the symmetric groups. The weak order is a partial order on the symmetric group. There is a natural linear basis of each symmetric group algebra,
related to the group basis by M\"obius inversion for the weak order. 
We describe the operad structure on this second basis: the surprising result is that each operadic composition is
a sum over an interval of the weak order.  We deduce that the coradical filtration is an operad filtration.
The  Lie operad, a suboperad of the associative operad, sits in the first component of the filtration. As a corollary to our results, we derive a simple
explicit expression for Dynkin's idempotent in terms of the second basis.

There are combinatorial procedures for constructing a planar binary tree from a
permutation, and a composition from a planar binary tree. These define set-theoretic
quotients of each symmetric group algebra. We show that they are non-symmetric operad quotients of the
associative operad. Moreover, the Hopf kernels of these quotient maps are non-symmetric suboperads
of the associative operad.
\end{abstract}

\maketitle

\section*{Introduction} \label{S:int}

One of the simplest symmetric operads is the associative operad
$\Ass$. This is an algebraic structure carried by the sequence of vector spaces $\Ass_n=\field S_n$, $n\geq 1$, where $S_n$ is the symmetric group on $n$ letters. In particular
this entails structure maps, for each $n,m\geq 1$ and $1\leq i\leq n$,
\[\Ass_n\otimes\Ass_m \map{\circ _i} \Ass_{n+m-1}\]
satisfying certain axioms (Section~\ref{S:associative}).
The Lie operad $\Lie$ is a symmetric suboperad of $\Ass$.

The space $\Ass_*:=\bigoplus_{n\geq 1}\field S_n$ carries the structure of a (non-unital) graded Hopf algebra, first defined by Malvenuto and Reutenauer~\cite{MR},
and studied recently in a number of works, including~\cite{AS,DHT01,LR2}.
It is known that $\Lie_*:=\bigoplus_{n\geq 1}\Lie_n$ sits inside the subspace of $\Ass_*$ consisting of
primitive elements for this Hopf algebra. This led us to consider whether this subspace is itself a suboperad of $\Ass$. In the process to
answering this question, we found a number of interesting results linking
the non-symmetric operad structure of $\Ass$ to the combinatorics of the symmetric groups,
and in particular to a partial order on $S_n$ known as the left weak
Bruhat order (or weak order, for simplicity).

Let $F_\sigma$ denote the standard basis element of $\field S_n$ corresponding to $\sigma\in S_n$. Define a new basis of $\field S_n$
by means of the formula
\[M_\sigma= \sum_{\sigma\leq \tau} \mu(\sigma,\tau) F_\tau\,,\]
where $\mu$ is the M\"obius function of the weak order (Section~\ref{S:weak}).
This basis was used in~\cite{AS} to provide a simple explicit
description of the primitive elements and the coradical filtration of the
Hopf algebra of Malvenuto and Reutenauer (as well as the rest of the
Hopf algebra structure). 

One of our main results, Theorem~\ref{T:operad-M}, provides an explicit
description for the non-symmetric operad structure of $\Ass$ (the maps $\circ_i$) in this basis.
We find that each $M_\sigma\circ_i M_\tau$ is a sum over the elements
of an interval in the weak order, which we describe explicitly. This result, and the combinatorics needed for its proof, are given in Sections~\ref{S:weak},~\ref{S:construction}, and~\ref{S:proof}.

Together with the results of~\cite{AS}, Theorem~\ref{T:operad-M} allows us to
conclude that the space of primitive elements is a non-symmetric suboperad of $\Ass$, and moreover that the coradical filtration is an operadic filtration.
These notions are reviewed in Section~\ref{S:coradical} and the result is obtained in Theorem~\ref{T:coradical}.

There are combinatorial procedures for constructing a planar binary tree with $n$ internal vertices from a
permutation in $S_n$, and a subset of $[n-1]$ from such a tree. 
The composite procedure associates to a permutation the set of its descents.
The behavior of these constructions with respect to the Hopf algebra structure is well-understood~\cite{GKal,LR2,Reu}.
 In Sections~\ref{S:descent} and~\ref{S:trees} we show that they lead
 to non-symmetric operad quotients of $\Ass$. These quotients possess
 partial orders and linear bases analogous to those of $\field S_n$,
 and the non-symmetric operad structure maps on the basis $M$ are
 again given by sums over intervals in the 
weak order (which degenerate to a point in the case of subsets). These results are
 obtained in Propositions~\ref{P:operad-M-Q} and~\ref{P:operad-M-Y}.
 Special attention is granted to the quotient of $\Ass$ defined by passing to
 descents. The Hopf kernel of this map is shown to be a non-symmetric suboperad of $\Ass$ in Proposition~\ref{P:hopfkernel}.

In Section~\ref{S:Lie-Dynkin} we turn to the symmetric operad $\Lie$. 
The subspace $\Lie_n$ of $\Ass_n$ is generated as right $S_n$-module by
a special element $\theta_n$ called Dynkin's idempotent. 
Our main result
here is a surprisingly simple expression for this element in the basis $M$ (Theorem~\ref{T:Dynkin}):
\[\theta_n=\sum_{\sigma\in S_n,\,\sigma(1)=1}M_\sigma\,.\]
We obtain this result by noting that the classical definition of $\theta_n$ can
be recast as the iteration of a certain operation
$\ebracket:\Ass\times\Ass\to\Ass$ that preserves the suboperad $\Lie$
and the non-symmetric suboperad of primitive elements, and by calculating an explicit expression for $\bracket{M_\sigma,M_\tau}$ (Proposition~\ref{P:bracket-M}).

\subsection*{Notation} The set $\{1,2,\ldots,n\}$ is denoted $[n]$.
We work over a commutative ring $\field$ of arbitrary characteristic. 
We refer to $\field$-modules as ``spaces''.

The symmetric group $S_n$ is the group of bijections $\sigma:[n]\to[n]$.
We denote permutations by the list of their values. Thus,
$\sigma=(\sigma_1,\ldots,\sigma_n)$ denotes a permutation $\sigma\in S_n$
whose value on $i$ is $\sigma_i$.

\section{The associative operad}\label{S:associative}

\subsection{Symmetric and non-symmetric operads}\label{S:non-symm}

A {\em non-symmetric unital operad} is a sequence of spaces $\calP=\{\calP_n\}_{n\geq 1}$  together with  linear maps
\[\calP_n\otimes\calP_m \map{\circ _i} \calP_{n+m-1}\,,\]
one for each $n,m\geq 1$ and $1\leq i\leq n$, and a distinguished element $1\in\calP_1$, such that
\begin{gather}
\begin{split}\label{E:operad-assoc}
(x\circ_i y)\circ _{j+m-1} z = (x\circ_j z)\circ_i y & \text{ \ if $1\leq
i<j\leq n$,} \\
(x\circ_i y)\circ _{i+j-1} z = x\circ_i (y\circ_j z)& \text{ \ if $1\leq 
i\leq n$ and $1\leq j\leq m$,}
\end{split}\\
\begin{split}\label{E:operad-unit}
x\circ_i 1 =  x  & \text{ \ if $1\leq i\leq n$,}\\
1\circ_1 x = x &
\end{split}
\end{gather}
for every  $x\in\calP_n$, $y\in\calP_m$, and $z\in\calP_l$. 
For various equivalent and related definitions, see~\cite{Lo96} or~\cite[II.1]{MSS}.



Given  permutations $\sigma=(\sigma_1,\ldots,\sigma_n)\in S_n$ and $\tau=(\tau_1,\ldots,\tau_m)\in S_m$, define a permutation $B_i(\sigma,\tau)\in S_{n+m-1}$ by
\begin{equation}\label{E:def-B}
B_i(\sigma,\tau):=(a_1,\ldots,a_{i-1},b_1,\ldots,b_m,a_{i+1},\ldots,a_n)
\end{equation}
where
\begin{equation}\label{E:def-B-2}
a_j:=\begin{cases}
\sigma_j & \text{ if }\sigma_j<\sigma_i \\
 \sigma_j+m-1        & \text{ if } \sigma_j>\sigma_i
\end{cases} \text{ \ \ and \ \ }
b_k:=\tau_k+\sigma_i-1\,.
\end{equation}
Note that $a_j\in [1,\sigma_i-1]\cup [\sigma_i+m,n+m-1]$ and $b_k\in
[\sigma_i,\sigma_i+m-1]$, so $B_i(\sigma,\tau)$ is indeed a
permutation. For instance, if $\sigma=(\blue{2},\blue{3},\blue{1},\blue{4})$, 
$\tau=(\red{\bf 2},\red{\bf 3},\red{\bf 1})$, and $i=2$ then 
$B_2(\sigma,\tau)=(\blue{2},\red{\bf 4},\red{\bf 5},\red{\bf 3},\blue{1},\blue{6})$. One may understand this construction in
terms of permutation matrices. Associate to $\sigma\in S_n$ the $n\times n$-matrix whose $(i,j)$-entry is Kronecker's $\delta_{i,\sigma_j}$. Then the matrix of $B_i(\sigma,\tau)$ is obtained by inserting
the matrix of $\tau$ in the $(\sigma_i,i)$ entry of the matrix of $\sigma$. In the above example,
\[\sigma=\mat{\blue{0} & \blue{0} & \blue{1} & \blue{0}\\
\blue{1} & \blue{0} & \blue{0} & \blue{0}\\
\blue{0} & \blue{1} & \blue{0} & \blue{0}\\
\blue{0} & \blue{0} & \blue{0} & \blue{1} }\,, \quad
\tau=\mat{\red{\bf 0} & \red{\bf 0} & \red{\bf 1} \\
\red{\bf 1} & \red{\bf 0} & \red{\bf 0}\\
\red{\bf 0} & \red{\bf 1} & \red{\bf 0} }\,, \text{ \ and \ }
B_2(\sigma,\tau)=\mat{\blue{0} & 0 & 0 & 0 & \blue{1} & \blue{0}\\
\blue{1} & 0 & 0 & 0 & \blue{0} & \blue{0}\\
0 & \red{\bf 0} & \red{\bf 0} & \red{\bf 1} & 0 & 0\\
0 & \red{\bf 1} & \red{\bf 0} & \red{\bf 0} & 0 & 0\\
0 & \red{\bf 0} & \red{\bf 1} & \red{\bf 0} & 0 & 0\\
\blue{0} & 0 & 0 & 0 & \blue{0} & \blue{1} }
\,.\]

A {\em symmetric unital operad} is a sequence of spaces $\calP$ with the
same structure as above, plus a right linear action of the symmetric group $S_n$ on $\calP_n$ such that
\begin{equation}\label{E:operad-equivariance}
(x\cdot \sigma)\circ_{i} (y\cdot\tau)=(x\circ_{\sigma(i)} y)\cdot
B_i(\sigma,\tau)
\end{equation}
for every $x\in \calP_n$, $y\in \calP_m$, $\sigma\in S_n$, $\tau\in S_m$,
and $1\leq i\leq n$.

Associated to any (symmetric or non-symmetric) operad  $\calP$, there
is the graded vector space
\[
\calP_*:=\bigoplus _{n\geq 1} \calP_n.
\]

There is a pair of adjoint functors 
\begin{center}
\{symmetric operads\} \raisebox{4pt}{$\map{\textit{\phantom{xati}forgetful\phantom{onx}}}$}\hspace*{-67pt}
\raisebox{-9pt}{$\overset{\leftmap{\phantom{regularization}}}{\scriptstyle{symmetrization}}$}     \{non-symmetric operads\}\,,
\end{center}
the symmetrization functor $\calS$ being left adjoint to the forgetful
functor $\calF$. Given a symmetric operad $\calP$, the non-symmetric
operad $\calF\calP$ is obtained
by forgetting the symmetric group actions. Conversely,
every non-symmetric operad $\calP$ gives rise to a symmetric operad $\calS\calP$
with spaces $\calS\calP_n:=\calP_n\otimes\field S_n$. These spaces are
equipped with the action of $S_n$ by right multiplication on the
second tensor factor. There is then a unique way to extend the
structure maps $\circ_i$ from 
$\calP$ to $\calS\calP$ in a  way that is compatible with the action. 


\medskip

An \emph{algebra} over a non-symmetric operad $\calP$ is a space $A$
together with structure maps 
\[
\calP_n\otimes A^{\otimes n} \to A
\]
subject to certain associativity and unitality conditions. An algebra over
a symmetric operad $\calP$ is a space $A$ as above for which the
structure maps factor through the quotient $\calP_n\otimes_{\field S_n} A^{\otimes n}$.

An algebra over a non-symmetric operad $\calP$ is the same thing as
an algebra over the symmetrization $\calS\calP$. On the other hand,
if $\calP$ is a symmetric operad, 
algebras over $\calP$ and algebras over $\calF\calP$ differ.

\subsection{The associative operad}\label{S:assoc}

Let $\Ass_n:=\field S_n$ be the group algebra of the symmetric group. The basis
element corresponding to a permutation $\sigma\in S_n$ is denoted $F_\sigma$.
This enables us to distinguish between various linear bases of $\Ass_n$, all of which
are indexed by permutations. In particular, a basis $M_\sigma$ is introduced in
Section~\ref{S:weak} below. 

The collection $\Ass:=\{\Ass_n\}_{n\geq 1}$ carries a structure of symmetric operad, uniquely determined by the requirements
\begin{align*}
F_{1_n}\circ_i F_{1_m}  =F_{1_{n+m-1}} &
\text{  \ for any $n,m\geq 1$, $1\leq i\leq n$,} \\
F_\sigma  =F_{1_n}\cdot\sigma &
\text{  \ for any $n\geq 1$, $\sigma\in S_n$,}
\end{align*}
where $1_n=(1,2,\ldots,n)$ denotes the identity permutation in $S_n$.
 In view of~\eqref{E:operad-equivariance},
the structure maps $\Ass_n\otimes\Ass_m \map{\circ _i} \Ass_{n+m-1}$
are given by
\begin{equation}\label{E:operad-F}
F_\sigma\circ _i F_\tau := F_{B_i(\sigma,\tau)}\,.
\end{equation}

This is the {\em associative} operad $\Ass$. It is a symmetric operad.
The non-symmetric operad $\calF\Ass$ is denoted $\Anon$.
We refer to $\Anon$ as the non-symmetric associative operad.
Algebras over $\Ass$ are associative algebras; we do not have
a simple description for the algebras over $\Anon$.

The {\em commutative} operad $\calC$ is the sequence of $1$-dimensional spaces $\field\{x_n\}$ with structure maps
\[x_n\circ_i x_m=x_{n+m-1}\,.\]
It is a symmetric operad with the trivial symmetric group actions.
The symmetric operad $\calS\calF\calC$ is the associative operad $\Ass$.

The Lie operad $\Lie$ is defined in Section~\ref{S:Lie}.
It is a symmetric suboperad of $\Ass$ and
algebras over $\Lie$ are Lie algebras. 
The non-symmetric operad $\calF\Lie$ is denoted $\Lnon$. It is a (non-symmetric) suboperad of $\Anon$.

\subsection{The weak order and the monomial basis}\label{S:weak}

The set of inversions of a permutation $\sigma\in S_n$ is
\[ \Inv(\sigma)\ :=\ \{(i,j)\in [n]\times[n] \mid i<j \text{ and }\sigma_i>\sigma_j\}\,.\]
The set of inversions determines the permutation. 

Let $\sigma,\tau\in S_n$. The {\em left weak Bruhat order} on $S_n$ is defined by
\[\sigma\leq\tau \iff \Inv(\sigma)\subseteq\Inv(\tau)\,.\]
This is a partial order on $S_n$. We refer to it as the {\em weak order} for simplicity.
For the Hasse diagram of the weak order on $S_4$, see~\cite[Figure 1]{AS}
or Figure~\ref{F:fibers} in Section~\ref{S:proof}.

Let $\sigma\leq\tau$ in $S_n$. The {\em M\"obius function} $\mu(\sigma,\tau)$ is defined by the recursion
\[\sum_{\sigma\leq\rho\leq\tau}\mu(\sigma,\rho)=
\begin{cases}
1 & \text{ if }\sigma=\tau\,, \\
 0        & \text{ if }\sigma<\tau\,.
\end{cases}\]
The M\"obius function of the weak order takes values in $\{-1,0,1\}$. 
Explicit descriptions  can be found in~\cite[Corollary 3]{Bj} 
or~\cite[Theorem 1.2]{Ede}.  
We will not need these descriptions.

The \emph{monomial} basis $\{M_\sigma\}$ of $\Anon_n$
is defined as follows~\cite[Section 1.3]{AS}. 
For each $n\geq 1$ and $\sigma\in S_n$, let

 \begin{equation}\label{E:def-monomial}
   M_\sigma := \sum_{\sigma\leq \tau} \mu(\sigma,\tau) F_\tau\,.
 \end{equation}
 For instance,
\[  M_{(4,1,2,3)}= F_{(4,1,2,3)}- F_{(4,1,3,2)}- F_{(4,2,1,3)}+F_{(4,3,2,1)}\,. \]

By M\"obius inversion, 
 \begin{equation} \label{E:fun-mon}
   F_\sigma = \sum_{\sigma\leq \tau} M_\tau\,.
 \end{equation}

Before giving the full description of the operad structure of $\Anon$ on the basis $\{M_\sigma\}$, consider one particular example. Using~\eqref{E:operad-F} and~\eqref{E:def-monomial}, one finds
by direct calculation that
\begin{equation}\label{E:ex-operad-M}
M_{(1,2,3)}\circ_2 M_{(2,1)}=M_{(1,3,2,4)}+M_{(1,4,2,3)}+M_{(2,3,1,4)}+M_{(2,4,1,3)}+M_{(3,4,1,2)}\,.\end{equation}

The five permutations appearing on the right hand side form an interval in
the weak order on $S_4$; the bottom element is $(1,3,2,4)$ and the top element is
$(3,4,1,2)$. This is a general fact. Fix $n,m\geq 1$ and $i\in[n]$. Below we define a map $T_i:S_n\times S_m\to S_{n+m-1}$ and we prove:

\begin{theorem}\label{T:operad-M}
For any $\sigma\in S_n$, $\tau\in S_m$, and $1\leq i\leq n$, 
\begin{equation}\label{E:operad-M}
M_\sigma\circ _i M_\tau = \sum_{B_i(\sigma,\tau)\leq\rho\leq T_i(\sigma,\tau)} M_{\rho}\,.
\end{equation}
\end{theorem}

A main ingredient in the proof of Theorem~\ref{T:operad-M} is the
construction of a map $P_i:S_{n+m-1}\to S_n\times S_m$ which is related to $B_i$ and $T_i$ through the following result.
We put on $S_n\times S_m$ the partial order obtained by taking the Cartesian product of the
weak orders on $S_n$ and $S_m$.

\begin{proposition}\label{P:connection} The maps 
\[P_i:S_{n+m-1}\to S_n\times S_m\,, \quad
B_i: S_n\times S_m\to S_{n+m-1}\,, \quad \text{ and } \quad
T_i: S_n\times S_m\to S_{n+m-1}\]
satisfy the following properties:
\begin{itemize}
\item[(i)] $P_i$ and $B_i$ are order-preserving.
\item[(ii)] $P_i\circ B_i=\Id = P_i\circ T_i$.
\item[(iii)] $B_i\circ P_i\leq \Id \leq T_i\circ P_i$.
\end{itemize}
\end{proposition}

The constructions of $T_i$ and $P_i$, and the necessary combinatorics,
are given in Section~\ref{S:construction}. The proofs of Theorem~\ref{T:operad-M}
and Proposition~\ref{P:connection}
are given in Section~\ref{S:proof}.

\subsection{Constructions with permutations}\label{S:construction}

First we set some notation. 
Given a sequence of distinct integers $a=(a_1,\ldots,a_n)$, we use $\{a\}$ to denote the underlying set $\{a_1,\ldots,a_n\}$.
The {\em standardization} of $a$ is the unique permutation $\st(a)$ of $[n]$ such that
$\st(a)_i<\st(a)_j$ if and only if $a_i<a_j$ for every $i,j\in[n]$.
The sequence $a$ is determined by the set $\{a\}$ and the permutation $\st(a)$.

Given two sets of integers $A$ and $B$ and an integer $m$, we write $A<m$ to indicate that $x<m$  for
every $x\in A$, and $A<B$ if $A< y$ for every $y\in
B$. 

Fix $n,m\geq 1$ and $i\in[n]$. All constructions below depend on $n$, $m$, and $i$, even when not explicitly mentioned.

Given $\rho\in S_{n+m-1}$,  define three sequences $\calL_j(\rho)$, $j=1,2,3$, by
\[\calL_1(\rho):=(\rho_1,\ldots,\rho_{i-1})\,,\ 
\calL_2(\rho):=(\rho_i,\ldots,\rho_{i+m-1})\,, \text{ \  and \ }
\calL_3(\rho):=(\rho_{i+m},\ldots,\rho_{n+m-1})\,.\] 

We proceed to define a map $T_i:S_n\times S_m\to S_{n+m-1}$.
Let $\sigma\in S_n$ and $\tau\in S_m$. 
 
 If $m=1$ we set $T_i(\sigma,1):=\sigma$. 
 
 Assume $m>1$. Let $L_j:=\calL_j(B_i(\sigma,\tau))$ for $j=1,2,3$. We
 will define $T_i(\sigma,\tau)$ by specifying the three sequences
 $\calL_j(T_i(\sigma,\tau))$, $j=1,2,3$.

Let $\eta_i:=\sigma_i+m-1$. Define
 \begin{align*}
 k_\sigma & :=\max\bigl\{j\in[0,i-1] \ \mid\  [\sigma_i-j,\sigma_i]\subseteq 
\{L_1\}\cup\{\sigma_i\} \bigr\}\,,\\
l_\sigma &:=\max\bigl\{j\in[0,n-i] \ \mid\  [\eta_i,\eta_i+j]\subseteq 
\{\eta_i\}\cup\{L_3\} \bigr\}\,.
\end{align*}
 The numbers $k_\sigma$ and $l_\sigma$ depend only on $\sigma$ (and $m$ and $i$), but not on $\tau$. The interval  $[\sigma_i-k_\sigma,\sigma_i-1]$ consists of the elements of
 $\{L_1\}$ that immediately  precede $\sigma_i$,
 and $[\eta_i+1,\eta_i+l_\sigma]$ consists of
 the elements of $\{L_3\}$ that immediately  follow
 $\eta_i$.
 The remaining elements of $\{L_1\}$ fall into two classes:
 $A_{1,b}^{\sigma}$, which consists of elements smaller than $\sigma_i-k_\sigma$, and $A_{1,t}^{\sigma}$, which consists of elements bigger than $\sigma_i$.  Similarly, the remaining elements of $\{L_3\}$ fall into two classes:
 $A_{3,b}^{\sigma}$, which consists of elements smaller than $\eta_i$, and $A_{3,t}^{\sigma}$, which consists of elements bigger than $\eta_i+l_\sigma$.

The interval $[\sigma_i-k_\sigma, \eta_i+l_\sigma]$ is thus partitioned as follows:

\[ \begin{picture}(400,30)(0,0)
   \put(0,0){\line(1,0){405}} 
   \put(0,-5){\line(0,1){10}} \put(120,-5){\line(0,1){10}}\put(150,-5){\line(0,1){10}}
   \put(300,-5){\line(0,1){10}}\put(330,-5){\line(0,1){10}}\put(405,-5){\line(0,1){10}}
  \blue{\put(50,30){\vector(-1,0){50}} \put(70,30){\vector(1,0){50}}
   \put(215,30){\vector(-1,0){65}}\put(235,30){\vector(1,0){65}}
    \put(360,30){\vector(-1,0){30}}   \put(375,30){\vector(1,0){30}}}
    \put(-12,12){$\sigma_i{-}k_\sigma$}    \put(107,12){$\sigma_i{-}1$}
    \put(145,12){$\sigma_i$} \put(295,12){$\eta_i$}
    \put(318,12){$\eta_i{+}1$}  \put(393,12){$\eta_i{+}l_\sigma$}      
     \put(57,28){$k_\sigma$}      \put(220,28){$m$}      \put(366,28){$l_\sigma$} 
\end{picture}\]

The numbers on top indicate the cardinality of each interval.

It follows from~\eqref{E:def-B-2} that $A_{3,b}^{\sigma}$ must consist of elements smaller
 than $\sigma_i-k_\sigma$, and $A_{1,t}^{\sigma}$ must consist of elements
 bigger than $\eta_i+l_\sigma$. Thus, we have
\begin{equation}\label{E:B-decomposition}
\begin{split}
\{L_1\}=&A_{1,b}^{\sigma}\cup [\sigma_i-k_\sigma,\sigma_i-1]\cup A_{1,t}^{\sigma}\,, \\
\{L_2\}=&[\sigma_i,\eta_i]\,, \\
\{L_3\}=&A_{3,b}^{\sigma}\cup [\eta_i+1,\eta_i+l_\sigma]\cup A_{3,t}^{\sigma}\,,
\end{split}
\end{equation}
with $A_{1,b}^{\sigma}$, $A_{3,b}^{\sigma}<\sigma_i-k_\sigma$ and
$\eta_i+l_\sigma<A_{1,t}^{\sigma}$, $A_{3,t}^{\sigma}$.

Consider now the following partition of the interval $[\sigma_i-k_\sigma, \eta_i+l_\sigma]$:

\[ \begin{picture}(400,50)(-15,-15)
   \put(-15,0){\line(1,0){405}} 
   \put(-15,-5){\line(0,1){10}} \put(15,-5){\line(0,1){10}}\put(90,-5){\line(0,1){10}}
   \put(120,-5){\line(0,1){10}}\put(210,-5){\line(0,1){10}}\put(240,-5){\line(0,1){10}}
   \put(360,-5){\line(0,1){10}}\put(390,-5){\line(0,1){10}}
   \blue{\put(43,30){\vector(-1,0){28}}   \put(60,30){\vector(1,0){30}}
    \put(150,30){\vector(-1,0){30}}   \put(180,30){\vector(1,0){30}}
   \put(290,30){\vector(-1,0){50}}   \put(310,30){\vector(1,0){50}}}
    \put(-30,12){$\sigma_i{-}k_\sigma$}    \put(-4,-15){$\sigma_i{-}k_\sigma{+}1$}
    \put(65,12){$\sigma_i{-}k_\sigma{+}l_\sigma$} \put(95,-15){$\sigma_i{-}k_\sigma{+}l_\sigma{+}1$}
    \put(175,12){$\eta_i{-}k_\sigma{+}l_\sigma{-}1$}  \put(225,-15){$\eta_i{-}k_\sigma{+}l_\sigma$}
     \put(336,12){$\eta_i{+}l_\sigma{-1}$}  \put(380,-15){$\eta_i{+}l_\sigma$}
     \put(50,28){$l_\sigma$}      \put(153,28){$m{-}2$}      \put(297,28){$k_\sigma$}      
\end{picture}\]

\smallskip

We define $T_i(\sigma,\tau)$ as the unique permutation in $S_{n+m-1}$ with sequences $L'_j:=\calL_j(T_i(\sigma,\tau))$ determined by
\[\st(L'_j)=\st(L_j)\]
 for $j\in\{1,2,3\}$, and
\begin{equation}\label{E:T-decomposition}
\begin{split}
\{L'_1\}=&A_{1,b}^{\sigma}\cup [\eta_i+l_\sigma-k_\sigma,\eta_i+l_\sigma-1] \cup A_{1,t}^{\sigma}\,, \\
\{L'_2\}=&
\{\sigma_i-k_\sigma\}\cup[\sigma_i-k_\sigma+l_\sigma+1,\eta_i+l_\sigma-k_\sigma-1] \cup \{\eta_i+l_\sigma\}\,, \\
\{L'_3\}=&A_{3,b}^{\sigma}\cup  
[\sigma_i-k_\sigma+1,\sigma_i-k_\sigma+l_\sigma]\cup A_{3,t}^{\sigma}\,. 
\end{split}
\end{equation}

\smallskip

Here is an example. Let $\sigma=(5,8,2,4,6,1,7,3)$, $\tau=(2,4,3,1)$, and $i=5$.
We have 
\[B_5(\sigma,\tau)=(5,11,2,4,7,9,8,6,1,10,3)\,,\]
and
\[L_1=(5,11,2,4)\,,\quad L_2=(7,9,8,6)\,,  \text{ \ and \ }L_3=(1,10,3)\,.\]
Since $\sigma_5=6$ and $\eta_5=9$, we have 
\[\{L_1\}=\{2\}\cup[4,5]\cup\{11\}\,,\quad
\{L_2\}=[6,9]\,, \text{ \ and \ }
\{L_3\}=\{1,3\}\cup[10,10]\cup\emptyset\,.\]
Hence $k_\sigma=2$, $l_\sigma=1$, and
\[\{L'_1\}=\{2\}\cup[8,9]\cup\{11\}\,,\quad
\{L'_2\}=\{4\}\cup[6,7]\cup\{10\}\,, \text{ \ and \ }
\{L'_3\}=\{1,3\}\cup[5,5]\cup\emptyset\,.\]
Finally,
\[L'_1=(8,11,2,9)\,,\quad L'_2=(6,10,7,4)\,,  \text{ \ and \ }L'_3=(1,5,3)\,,\]
so
\[T_5(\sigma,\tau)=(8,11,2,9,6,10,7,4,1,5,3)\,.\]


\smallskip

Next we define a map $P_i:S_{n+m-1}\to S_n\times S_m$ (the map depends on $n$ and $m$ but this is omitted from the notation). Let  $\rho\in S_{n+m-1}$.
If $m=1$ we set $P_i(\rho):=(\rho,1)$. 

Assume $m\geq 2$. Define
\[u^\rho:=\min\{\calL_2(\rho)\} \text{ \ and \ }v^\rho:=\max\{\calL_2(\rho)\}\,.\]
Write 
\begin{equation}\label{E:rho-decomposition}
\{\calL_1(\rho)\}=C_{1,b}^\rho\cup C_{1,m}^\rho \cup C_{1,t}^\rho \text{ \ and \ }
\{\calL_3(\rho)\}=C_{3,b}^\rho\cup C_{3,m}^\rho \cup C_{3,t}^\rho
\end{equation}
 with 
 \[C_{j,b}^\rho<u^\rho< C_{j,m}^\rho
< v^\rho< C_{j,t}^\rho\]
for $j=1,3$.
 
Thus $[u^\rho,v^\rho]=C_{1,m}^\rho\cup C_{3,m}^\rho\cup \{\calL_2(\rho)\}$. If $n_{1,m}^\rho$ (resp. $n_{3,m}^\rho$) denotes
the number of elements in $C_{1,m}^\rho$ (resp. $C_{3,m}^\rho$), then $v^\rho-u^\rho+1=m+n_{1,m}^\rho+
n_{3,m}^\rho$.

Define a permutation $B\rho\in S_{n+m-1}$ by
$\st(\calL_j(B\rho))=\st(\calL_j(\rho))$ for $j=1,2,3$ and
\begin{equation}\label{E:BP-decomposition}
\begin{split}
\{\calL_1(B\rho)\}:=&C_{1,b}^\rho\cup [u^\rho,u^\rho+n_{1,m}^\rho-1]\cup C_{1,t}^\rho\,, \\
\{\calL_2(B\rho)\}:=&[u^\rho+n_{1,m}^\rho,v^\rho-n_{3,m}^\rho ]\,,\\
\{\calL_3(B\rho)\}:=&C_{3,b}^\rho\cup [v^\rho-n_{3,m}^\rho+1,v^\rho]\cup C_{3,t}^\rho\,.
\end{split} 
\end{equation}

Finally, define
\begin{equation}\label{E:def-P}
P_i(\rho):=(\st(\calL_1(B\rho),u^\rho+n_{1,m}^\rho,\calL_3(B\rho)),\st(\calL_2(\rho)))\,.
\end{equation}
By construction,
\begin{equation}\label{E:rel-B-P}
B_i(P_i(\rho))=B\rho\,.
\end{equation}

\subsection{Proofs of Theorems~\ref{T:operad-M}
and Proposition~\ref{P:connection}}\label{S:proof}

\begin{proof}[Proof of Proposition~\ref{P:connection}]
The case $m=1$ is trivial. Assume $m>1$.

Let $I_1:=[1,i-1]$, $I_2:=[i,i+m-1]$,
$I_3:=[i+m,n+m-1]$. For $\alpha\in I_1$, $\beta\in I_2$, and $\gamma\in I_3$, let $\tilde\alpha:=\alpha$,
$\tilde\beta:=\beta-i+1$, and $\tilde \gamma:=\gamma-m+1$. Choose
   $\alpha_j,\beta_j\in I_j$ for each $j=1,2,3$.

The following assertions are easy to check: for $k,l\not=2$ and  $k\leq l$,
\begin{align*}
(\alpha_k,\beta_l)\in \Inv(B_i(\sigma,\tau))& \Longleftrightarrow
(\tilde\alpha_k,\tilde\beta_l)\in \Inv(\sigma)\,, \\
 (\alpha_2,\beta_2)\in \Inv(B_i(\sigma,\tau))& \Longleftrightarrow
(\tilde\alpha_2,\tilde\beta_2)\in \Inv(\tau)\,,\\
(\alpha_1,\beta_2)\in  \Inv(B_i(\sigma,\tau))& \Longleftrightarrow
(\alpha_1,i)\in\Inv(\sigma)\,,\\
(\alpha_2,\beta_3)\in  \Inv(B_i(\sigma,\tau)) &\Longleftrightarrow
(i,\tilde\beta_3)\in\Inv(\sigma)\,.
\end{align*}
It follows that
\begin{equation}\label{E:B-order-preserving}
B_i(\sigma,\tau)\leq
B_i(\sigma',\tau')\Leftrightarrow (\sigma,\tau)\leq(\sigma',\tau')\,.
\end{equation}
Hence,  $B_i$ is injective and order-preserving. 

Given $\rho$  in $S_{n+m-1}$ and $h\in\{1,3\}$, let   
\[I_\rho^h:=\{k_h\in I_h \mid  \exists\ i_2,j_2\in I_2, 
\rho(i_2)<\rho(k_h)<\rho(j_2)\}= \{k_h\in I_h \mid \rho(k_h)\in C_{h,m}^\rho\}\]
and 
\[J_\rho:=\Inv(\rho)\cap (I_\rho^1\times I_2\cup 
I_2\times I_\rho^3\cup I_\rho^1\times I_\rho^3).\]
By~\eqref{E:rho-decomposition} and~\eqref{E:BP-decomposition} we have
\begin{equation}\label{E:inversionJ}
\Inv(B_iP_i(\rho))=\Inv(\rho)-J_\rho\,,
\end{equation}
which implies
\[ B_i\circ P_i \leq \Id\,.\]

Let $\rho'\in S_{n+m-1}$ be  such that $\rho\leq \rho'$. Let $(k,l)$ be in $J_{\rho'}\cap \Inv(\rho)$.
If $k\in I_{\rho'}^1$, there exist $i_2,j_2\in I_2$ such
that
$\rho'(i_2)<\rho'(k)<\rho'(j_2)$. Hence,
$(k,j_2)\not\in\Inv(\rho')$, so $(k,j_2)\not\in\Inv(\rho)$. But
$(k,l)\in \Inv(\rho)$ implies $\rho(l)<\rho(k)<\rho(j_2)$ hence
$k\in I_\rho^1$ if $l\in I_2$. If $l\in I_{\rho'}^3$, there exists $j\in
I_2$ such that $\rho'(j)<\rho'(l)$, hence
$\rho(j)<\rho(l)<\rho(k)<\rho(j_2)$ and $k\in I_\rho^1$. The same argument applies to $l$:  if $l\in
I_{\rho'}^3$ then $l\in I_\rho^3$ when $k\in I_\rho^1\cup I_2$. As a consequence,  $J_{\rho'}\cap \Inv(\rho)\subseteq
J_\rho$
and $\Inv(B_iP_i(\rho))\subseteq\Inv(B_iP_i(\rho'))$.
By (\ref{E:B-order-preserving}) we get that $P_i$ is
order-preserving.

If $\rho=B_i(\sigma,\tau)$ then $I_\rho^1=I_\rho^3=\emptyset$ and
$\Inv(B_iP_i(\rho))=\Inv(\rho)$. It follows that
$B_iP_i(\rho)=B_i(\sigma,\tau)$, and the injectivity of $B_i$ implies
\[ P_i\circ B_i=\Id\,.\]

We  compare $B_i\circ P_i\circ T_i$ to $B_i$. Given $(\sigma,\tau)\in S_n\times S_m$, 
it suffices to compare $B\rho$ to $B_i(\sigma,\tau)$
with $\rho:=T_i(\sigma,\tau)$, in view of~\eqref{E:rel-B-P}. Since the
standardization of the lists $\calL_j$ for these two permutations are the same, it suffices 
to compare the sets $\{\calL_j(B\rho)\}$ and
$\{\calL_j(B_i(\sigma,\tau)\}$. These sets coincide in view of
\eqref{E:B-decomposition}, \eqref{E:T-decomposition}, \eqref{E:rho-decomposition},  and
\eqref{E:BP-decomposition}. Since $B_i$ is injective, it follows that
\[P_i\circ T_i=\Id\,.\]

With the notation of \eqref{E:B-decomposition}, let 
\begin{align*}
I_1^\sigma & =\{\alpha\in I_1 \mid \sigma_i-k_\sigma\leq 
B_i(\sigma,\tau)(\alpha)\leq \sigma_i-1\}\,,\\ 
I_3^\sigma & =\{\beta\in I_3 \mid \eta_i+1\leq B_i(\sigma,\tau)(\beta)\leq \eta_i+l_\sigma\}\,,\\
K_{B_i(\sigma,\tau)} & =\{(\alpha,\beta)\in I_1^\sigma\times I_2 \mid B_i(\sigma,\tau)(\beta)<\eta_i\}\\
&\  \cup \{(\alpha,\beta)\in I_2\times I_3^\sigma \mid \sigma_i<B_i(\sigma,\tau)(\alpha)\}\cup 
I_1^\sigma\times I_3^\sigma\,.
\end{align*}

The inversion set of $T_i(\sigma,\tau)$ is
\begin{equation}\label{E:inversionK}
\Inv(T_i(\sigma,\tau))=\Inv(B_i(\sigma,\tau))\cup K_{B_i(\sigma,\tau)}\,.
\end{equation}
In view of (\ref{E:inversionJ}) and (\ref{E:inversionK}),  
to show that $T_i\circ P_i\geq \Id$ it suffices to prove that
$J_\rho\subseteq K_{B_iP_i(\rho)}$.
Recall from \eqref{E:rel-B-P} that $B_iP_i(\rho)=B\rho$. Let $(\sigma,\tau):=P_i(\rho)$. Then $k_\sigma\geq
n_{1,m}^\rho$ and $l_\sigma\geq n_{3,m}^\rho$. These imply that if 
$j\in I^\alpha_\rho$ then $j\in I_\alpha^\sigma$. Assume $(j,l)\in I_\rho^1\times I_2$ is 
an inversion for $\rho$. Then $\rho(l)<\rho(j)<v^{\rho}$ and $B\rho(l)<v^{\rho}-n_{3,m}^\rho=\sigma_i+m-1=\eta_i$.
It follows that $(j,l)$ lies in $K_{B\rho}$. The same argument holds if $(j,l)\in I_2\times I_\rho^3$.
Hence $J_\rho\subseteq K_{B_iP_i(\rho)}$ and
\[T_i\circ P_i\geq \Id.\]
\end{proof}

\begin{figure}[hb]
  $$  \epsfxsize=3in\epsfbox{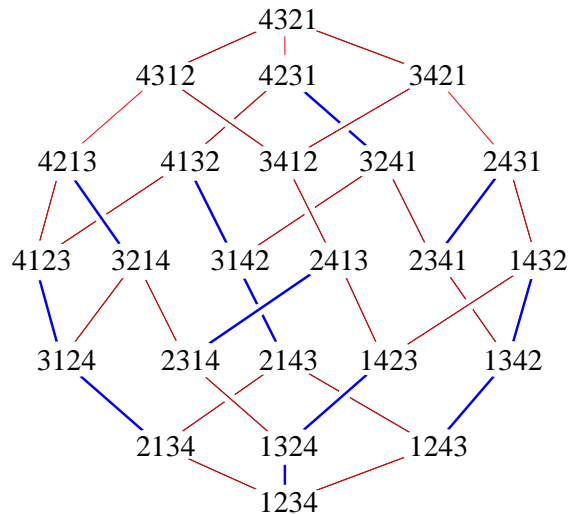}$$
  \caption{The fibers of $P_i:S_4\to S_3\times S_2$}\label{F:fibers}
\end{figure}

It follows from Proposition~\ref{P:connection} that the fiber of the map
$P_i$ over each pair $(\sigma,\tau)$ is an interval in the poset $S_{n+m-1}$, with bottom
element $B_i(\sigma,\tau)$ and top element $T_i(\sigma,\tau)$. 
Note that this data does {\em not} define a {\em lattice congruence} in the
sense of~\cite{Rea}, since the map $T_i$ is not order-preserving. The fact that $B_i$ is order-preserving 
has the following additional consequence: if $x\leq y$ in $S_{n+m-1}$, then the bottom element in the fiber of $P_i$  which contains $x$  is less than or equal to the bottom element in the fiber of $P_i$ which contains $y$. 
The fibers of $P_i$ are shown in Figure~\ref{F:fibers}, for $n=3$, $m=2$, $i=1$
(elements joined by a thick edge belong to the same fiber).

Equation~\eqref{E:operad-M} may thus be restated as follows:
\begin{equation}\label{E:operad-M-P}
M_\sigma\circ _i M_\tau = \sum_{P_i(\rho)=(\sigma,\tau)} M_{\rho}\,.
\end{equation}

\begin{proof}[Proof of Theorem~\ref{T:operad-M}]
Define a map $\tilde{\circ}_i$ by
\[M_\sigma\tilde{\circ}_i M_\tau:= \sum_{P_i(\rho)=(\sigma,\tau)} M_{\rho}\,.\]
Then, by~\eqref{E:fun-mon},
\[F_\sigma\tilde{\circ}_i F_\tau =\sum_{\sigma\leq\sigma' \atop \tau\leq\tau'} M_{\sigma'}\tilde{\circ}_i M_{\tau'}=
\sum_{(\sigma,\tau)\leq(\sigma',\tau') \atop P_i(\rho)=(\sigma',\tau')} M_{\rho}=
\sum_{(\sigma,\tau)\leq P_i(\rho)} M_\rho \,.
\]

 Proposition~\ref{P:connection} implies that $B_i$ is left adjoint to $P_i$, that is
\[B_i(\sigma,\tau)\leq \rho \Longleftrightarrow (\sigma,\tau)\leq P_i(\rho).\]
Therefore, 
\[F_\sigma\tilde{\circ}_i F_\tau = \sum_{B_i(\sigma,\tau)\leq\rho} M_{\rho}=F_{B_i(\sigma,\tau)}\,.\]
Thus $F_\sigma\tilde{\circ}_i F_\tau =F_\sigma\circ_i F_\tau$ and hence also, by linearity, $M_\sigma\tilde{\circ}_i M_\tau =M_\sigma\circ_i M_\tau$, which proves~\eqref{E:operad-M-P} and~\eqref{E:operad-M}.
\end{proof}

\section{A filtration of the non-symmetric associative operad}\label{S:coradical}

The space $H:=\bigoplus_{n\geq 1}\field S_n$ carries a graded Hopf algebra structure,
first introduced by Malvenuto and Reutenauer~\cite{MR}. The component of degree $n$ is $\field S_n$. We are interested in the graded coalgebra structure, which is defined for $\sigma\in S_n$ by
\begin{equation}\label{E:coproduct-F}
\Delta(F_{\sigma})=\sum_{i=1}^{n-1} F_{\st(\sigma_1,\,\ldots,\,\sigma_i)}\otimes
                                   F_{\st(\sigma_{i+1},\,\ldots,\,\sigma_n)}\,.
\end{equation}
(See \ref{S:construction} for the notion of standardization.)
This structure is studied in various recent works, including~\cite{AS,DHT01,LR2,MR}  (these references deal with the counital version of this coalgebra, which is obtained by adding a copy of the base ring in degree $0$ and adding the terms $1\otimes F_\sigma$ and $F_\sigma\otimes 1$ to~\eqref{E:coproduct-F}). 

We use $\sigma^i_{(1)}$ and $\sigma^i_{(2)}$ to
 denote the permutations $\st(\sigma_1,\,\ldots,\,\sigma_i)$ and $\st(\sigma_{i+1},\,\ldots,\,\sigma_n)$, respectively.
In particular, $\sigma^0_{(2)}=\sigma$ and 
$\sigma^n_{(1)}=\sigma$.

The iterated coproducts are defined by
\[\Delta^{(1)}:=\Delta \text{ \ and \ } \Delta^{(k+1)}:=(\Delta\otimes\id^{\otimes k})\circ\Delta^{(k)}\,.\]
Let $k\geq 1$. The $k$-th component of the {\em coradical filtration} is
\[H^{(k)}:=\ker(\Delta^{(k)}:H\to H^{\otimes (k+1)})\,.\]
The first component $H^{(1)}$ is the space of primitive elements of $H$. Note that
$H^{(1)}\subseteq H^{(2)}\subseteq H^{(3)}\subseteq\cdots$.

Since $\Delta$ is degree-preserving, each space $H^{(k)}$ is graded by setting
$H^{(k)}_n:=H^{(k)}\cap\field S_n$. For $k\geq 0$, let $\Anon^{(k)}$ denote  the sequence of spaces $\{H^{(k+1)}_n\}_{n\geq 1}$. It makes sense to wonder if  the sequence $\{\Anon^{(k)}\}_{k\geq 0}$
is a filtration of the non-symmetric associative operad $\Anon$. This question is motivated by the well-known fact that 
$\Anon^{(0)}$ contains the non-symmetric Lie operad $\Lnon$, a  suboperad of $\Anon$ (see Section~\ref{S:Lie} for more on this).

There does not seem to be a simple relation between the coalgebra structure of $H$ and
the operad structure of $\Anon$: while the following formula is easily verified
\[\Delta(F_\sigma\circ_i F_\tau)=\sum_{j=1}^{i-1} F_{\sigma^j_{(1)}}\otimes
(F_{\sigma^j_{(2)}}\circ_{i-j} F_{\tau})+
\sum_{k=1}^{m-1} (F_{\sigma^i_{(1)}}\circ_i F_{\tau^k_{(1)}}) \otimes
(F_{\sigma^{i-1}_{(2)}}\circ_{1} F_{\tau^k_{(2)}})+
\sum_{j=i}^{n-1} (F_{\sigma^j_{(1)}}\circ_{i} F_\tau) \otimes
F_{\sigma^j_{(2)}}\,, \]
its right hand side cannot be expressed entirely in terms of $\Delta(F_\sigma)$ and $\Delta(F_\tau)$, so it is not clear how this formula translates to other linear  bases of the group algebras. In spite of this, we have the following result. 

\begin{theorem}\label{T:coradical} The sequence $\{\Anon^{(k)}\}_{k\geq 0}$ is a filtration  of the non-symmetric associative operad $\Anon$, i.e.,
\begin{equation}\label{E:coradical}
\Anon^{(k)}_n \circ_i \Anon^{(h)}_m \subseteq \Anon^{(k+h)}_{n+m-1} \quad \text{  for any $n,m\geq 1$, $k,h\geq 0$, and $i=1,\ldots,n$.}
\end{equation}
In particular, the space of primitive elements $\Anon^{(0)}$ is a non-symmetric suboperad of $\Anon$.
\end{theorem}

We do not have a simple description of the algebras over the non-symmetric operad $\Anon^{(0)}$.

The proof of Theorem~\ref{T:coradical} makes use of an explicit
description of the coradical filtration in terms of the $M$-basis found in~\cite{AS}.  A permutation $\sigma\in S_n$ has a global descent at a position
 $p\in[n{-}1]$ if
  \[ \forall\ i\leq p \text{ and }j\geq p{+}1\,,\ \sigma_i>\sigma_j\,.\]
 Let $\GDes(\sigma)\subseteq[n{-}1]$ be the set of global descents of $\sigma$.
 Corollary 6.3 in~\cite{AS} states that the set
 \begin{equation}\label{E:k-global}
\{ M_\sigma \mid \sigma \text{ has at most $k$ global descents}\}
\end{equation}
 is a linear basis of $H^{(k+1)}$.

Given a set of integers $\setS$ and an integer $m$,  $\setS+m$ is the set with elements $s+m$ for $s\in \setS$.

\begin{lemma}\label{L:coradical} 
For any $\sigma\in S_n$, $\tau\in S_m$, and $i\in[n]$,
\begin{multline}\label{E:GDes-T1}
\GDes\bigl(T_i(\sigma,\tau)\bigr)\cap\bigl([1,i-1]\cup[i+m-1,n+m-2]\bigr)=\\
(\GDes(\sigma)\cap[1,i-1])\cup \bigl((\GDes(\sigma)\cap[i,n-1])+m-1\bigr),
\end{multline}
\begin{equation}\label{E:GDes-T}
\GDes\bigl(T_i(\sigma,\tau)\bigr)\cap [i,i+m-2]=\begin{cases}
\GDes(\tau)+i-1 & \text{ if } A_{1,b}^\sigma=\emptyset \text{ and } A_{3,t}^\sigma=\emptyset\,, \\
\emptyset        & \text{ if not.}
\end{cases}\\
\end{equation}
\end{lemma}
\begin{proof} Equation~\eqref{E:GDes-T1} follows easily from the definition of $T_i$. 
Let $\rho:=T_i(\sigma,\tau)$. 
Let $j\in [i,i+m-2]$ be a global descent of $\rho$. There exist $\alpha\in [i,j]$ and $\beta\in[j+1,i+m-1]$
such that $\rho(\alpha)=\eta_i+l_\sigma$ and   $\rho(\beta)=\sigma_i-k_\sigma$. Since $j$ is a global descent, one has
\[ \rho(u)>\sigma_i-k_\sigma \text{ \ and \ }
\rho(v)<\eta_i+l_\sigma\]
for all $u\in  [1,i-1]$ and $v\in [i+m,n+m-1]$. 
Hence $A_{1,b}^\sigma=\emptyset$, $A_{3,t}^\sigma=\emptyset$ and $j-i+1$ is a global descent of $\tau$. Conversely,
assume that  $A_{1,b}^\sigma=\emptyset$ and $A_{3,t}^\sigma=\emptyset$, and let $j$ be a global descent of $\tau$. 
There exist $\alpha\in [1,j]$ and $\beta\in[j+1,m]$
such that $\tau(\alpha)=m$ and   $\tau(\beta)=1$, so
$\rho(\alpha+i-1)=\eta_i+l_\sigma$ and  
$\rho(\beta+i-1)=\sigma_i-k_\sigma$. Then, using~\eqref{E:T-decomposition}, 
\[\forall\, s\leq j+i-1,\  \rho(s)>\rho(\beta+i-1)\implies
\rho(s)>\sigma_i-k_\sigma+l_\sigma 
\text{ and }\] 
\[\forall\, t>j+i-1,\ \rho(t)<\rho(\alpha+i-1)\implies \rho(t)<\eta_i-k_\sigma+l_\sigma\]
implies
\[\forall\, s\leq j+i-1,\ \rho(s)> \{\calL_3(\rho)\} \text{ and }  \forall t>j+i-1, \{\calL_1(\rho)\}>\rho(t).\]
Since $\{\calL_1(\rho)\}>\{\calL_3(\rho)\}$, \eqref{E:GDes-T} follows.
\end{proof}

\begin{proof}[Proof of Theorem~\ref{T:coradical}]
It follows from~\eqref{E:GDes-T} that the number of global descents
of $T_i(\sigma,\tau)$ is at most the sum of the numbers of global descents of $\sigma$ and $\tau$. Thus,  if $\sigma$ has at most $k$ global descents and
$\tau$ has at most $h$ global descents then $T_i(\sigma,\tau)$ has at most
$k+h$ global descents. Now, the map $\GDes$ is an order-preserving application between the weak order on $S_n$ and the poset of subsets of $[n-1]$ under inclusion~\cite[Proposition 2.13]{AS}. Therefore, the number of global descents of $\rho$ is at most $k+h$ for any permutation $\rho\leq T_i(\sigma,\tau)$. Hence, all permutations
appearing in the right hand side of~\eqref{E:operad-M} have at most $k+h$ global
descents, which proves~\eqref{E:coradical}, in view of~\eqref{E:k-global}.
\end{proof}

As already mentioned,  the set
$\{ M_\sigma \mid \sigma \text{ has at most $k$ global descents}\}$
 is a linear basis of $\Anon^{(k)}_*$.
In this sense, the operad $\Anon$ is filtered by the number of global descents.
Note that this is {\em not} an operad grading though.
Referring to the calculation of $M_{(1,2,3)}\circ_2 M_{(2,1)}$, we see that $(1,2,3)$ has $0$ global descents, $(2,1)$ has $1$, and the first four  permutations appearing in the right hand side of~\eqref{E:ex-operad-M} have $0$ global
descents, while the last one has $1$.

\section{Quotient operad: descent sets}\label{S:descent}

There is a combinatorial procedure for constructing a subset of $[n-1]$  from a
permutation of $[n]$. In this section we show that it leads to an
operad quotient of the non-symmetric associative operad $\Anon$. 

 A permutation $\sigma\in S_n$ has a  descent at a position
 $p\in[n{-}1]$ if
  \[ \sigma_p>\sigma_{p+1}\,.\]
 Let $\Des(\sigma)\subseteq[n{-}1]$ be the set of  descents of $\sigma$. 
 Note that $\GDes(\sigma)\subseteq\Des(\sigma)$. 
 
 For each $n\geq 1$, let $Q_n$ denote the poset of subsets of $[n-1]$ under inclusion (the Boolean poset). The map 
 $\Des:S_n\to Q_n$ is an order-preserving application from the weak order on $S_n$ to the Boolean poset $Q_n$. Let $\calQ$ denote the sequence of spaces $\{\field Q_n\}_{n\geq 1}$. The basis element of $\field Q_n$ corresponding to a subset $\setS\subseteq[n-1]$ is denoted $F_\setS$. 
 
 Given a set of integers $\setS$ and an integer $p$, let
 \[\setS+p:=\{s+p \mid s\in\setS\}\,.\]
 Given $\setS\subseteq [n-1]$, $\setT\subseteq [m-1]$, and $i\in[n]$, let
 \begin{equation}\label{E:def-B-Q}
B_i(\setS,\setT):=(\setS\cap\{1,\ldots,i-1\})\cup(\setT+i-1)\cup(\setS\cap\{i,\ldots,n-1\}+m-1)\subseteq [n+m-2]\,.
\end{equation}
 
 \begin{lemma}\label{L:descent} For any $\sigma\in S_n$, $\tau\in S_m$, and $i\in[n]$,
 \begin{equation}\label{E:descent}
\Des\bigl(B_i(\sigma,\tau)\bigr)=B_i\bigl(\Des(\sigma),\Des(\tau)\bigr)\,.
\end{equation}
\end{lemma}
\begin{proof} Write $B_i(\sigma,\tau)=(a_1,\ldots,a_{i-1},b_1,\ldots,b_m,a_{i+1},\ldots,a_n)$, as in~\eqref{E:def-B}. Since
\[\st(a_1,\ldots,a_{i-1})=\st(\sigma_1,\ldots,\sigma_{i-1}),\qquad \st(a_{i+1},\ldots,a_n)=\st(\sigma_{i+1},\ldots,\sigma_n)\]
\[\text{ \ and \ } \st(b_1,\ldots,b_m)=\st(\tau_1,\ldots,\tau_m),\]
the two sets appearing in~\eqref{E:descent} contain the same elements from 
\[ [1,i-2]\cup [i,i+m-2]\cup [i+m,n+m-2]\,.\]
If $i-1\in \Des\bigl(B_i(\sigma,\tau)\bigr)$ then $a_{i-1}>b_1$, that is, by \eqref{E:def-B-2}, $a_{i-1}>\tau_1+\sigma_i-1\geq\sigma_i$.
This can only occur if $\sigma_{i-1}>\sigma_i$, i.e., $i-1\in \Des(\sigma)$. Also,
  \begin{multline*}
  i+m-1\in \Des(B_i(\sigma,\tau)) \iff b_m=\tau_m+\sigma_i-1>a_{i+1} 
  \iff \sigma_i>\sigma_{i+1} \\
  \iff i\in\Des(\sigma)\iff i+m-1\in B_i\bigl(\Des(\sigma),\Des(\tau)\bigr)\,.
  \end{multline*}
  This completes the proof of~\eqref{E:descent}.
\end{proof}
 
 \begin{proposition}\label{P:operad-F-Q}
The sequence $\calQ$ is a non-symmetric operad under the structure maps
\begin{equation}\label{E:operad-F-Q}
F_\setS\circ _i F_\setT := F_{B_i(\setS,\setT)}\,.
\end{equation}
Moreover, the map $\sfD:\Anon\to\calQ$, $\sfD(F_\sigma):=F_{\Des(\sigma)}$, is
a surjective morphism of non-symmetric operads.
\end{proposition}
\begin{proof} According to~\eqref{E:operad-F} and Lemma~\ref{L:descent},
$\calQ$ inherits the structure of $\Anon$. Hence $\calQ$ is a non-symmetric operad and $\sfD$ a morphism.
\end{proof}
 
 In analogy with~\eqref{E:def-monomial}, the monomial basis of $\field Q_n$ is defined by  
 
 \begin{equation}\label{E:def-monomial-Q}
   M_\setS := \sum_{\setS\subseteq \setT} (-1)^{\#(T\setminus\setS)} F_\setT\,.
 \end{equation}
 For instance, with $n=4$,
\[  M_{\{1\}}= F_{\{1\}}- F_{\{1,2\}}- F_{\{1,3\}}+F_{\{1,2,3\}}\,. \]

By M\"obius inversion, 
 \begin{equation} \label{E:fun-mon-Q}
   F_\setS = \sum_{\setS\subseteq \setT} M_\setT\,.
 \end{equation}
 
 The operad structure of $\calQ$ takes the same form in the $M$-basis as in the $F$-basis.
 
 \begin{proposition}\label{P:operad-M-Q}
 For any $\setS\subseteq[n-1]$, $\setT\subseteq[m-1]$, and $i\in[n]$,
 \begin{equation}\label{E:operad-M-Q}
M_\setS\circ_i M_\setT=M_{B_i(\setS,\setT)}\,.
\end{equation}
\end{proposition}
\begin{proof} Define a map $\tilde{\circ}_i$ by
\[M_\setS\tilde{\circ}_i M_\setT:=M_{B_i(\setS,\setT)}\,.\]
Then, by~\eqref{E:fun-mon-Q},
\[F_\setS\tilde{\circ}_i F_\setT =\sum_{\setS\subseteq\setS' \atop \setT\subseteq\setT'} M_{\setS'}\tilde{\circ}_i M_{\setT'}=
\sum_{\setS\subseteq\setS' \atop \setT\subseteq\setT'}M_{B_i(\setS',\setT')}\,.
\]
Suppose one is given $B_i(\setS,\setT)$, $n$, $m$, and $i$. An inspection of~\eqref{E:def-B-Q} reveals that one can then determine $\setS$ and $\setT$.
Fix $\setS$ and $\setT$. It follows that the map $B_i$
is a bijection between the set 
\[\{(\setS',\setT') \mid \setS\subseteq\setS'\subseteq[n-1] \text{ and } \setT\subseteq\setT'\subseteq[m-1]\}\]
and the set
\[\{\setR \mid B_i(\setS,\setT)\subseteq \setR\subseteq[n+m-2]\}\,.\]
Therefore,
\[F_\setS\tilde{\circ}_i F_\setT =\sum_{B_i(\setS,\setT)\subseteq \setR} M_\setR=
F_{B_i(\setS,\setT)}\,.\]
Thus $F_\setS\tilde{\circ}_i F_\setT=F_\setS\circ_i F_\setT$ and hence also, by linearity,
$M_\setS\tilde{\circ}_i M_\setT=M_\setS\circ_i M_\setT$, which proves~\eqref{E:operad-M-Q}.
\end{proof}

\begin{corollary}\label{C:isom-Q} The map
$F_\setS\mapsto M_\setS$
defines an automorphism of the non-symmetric operad $\calQ$.
\end{corollary}
\begin{proof} This is an immediate consequence of~\eqref{E:operad-F-Q} and~\eqref{E:operad-M-Q}.
\end{proof}

The analogous formula to~\eqref{E:operad-M-Q} at the level of the operad $\Anon$ is~\eqref{E:operad-M}. The latter involves a sum over an interval, while in the former the interval has degenerated to a single point. For instance, 
the formula that corresponds to~\eqref{E:ex-operad-M} is simply
($n=3$, $m=2$)
\begin{equation}\label{E:ex-operad-M-Q}
M_{\emptyset}\circ_2 M_{\{1\}}=M_{\{2\}}\,.
\end{equation}

Note that $F_\sigma\mapsto M_\sigma$ does {\em not} define an automorphism of the operad $\Anon$, in contrast to Corollary~\ref{C:isom-Q}.
A further comparison between the two formulas leads to the following observations.

First, let us recall the expression of the map $\sfD$ on the $M$-bases of $\Anon$ and $\calQ$. We say that a permutation $\sigma\in S_n$ is closed if $\Des(\sigma)=\GDes(\sigma)$, as in~\cite[Definition 7.1]{AS}. The map $\Des$ (and also the map $\GDes$) defines an isomorphism between the subposet of $S_n$ consisting of closed permutations and the Boolean poset $Q_n$~\cite[Proposition 2.11]{AS}.

The map $\sfD$ is given as follows~\cite[Theorem 7.3]{AS}
\begin{equation}\label{E:D-monomial}
 \sfD(M_\sigma)\ =\ \begin{cases}
           M_{\Des(\sigma)} & \text{if $\sigma$ is closed,}\\
                     0 & \text{if not.}
        \end{cases}
\end{equation}

\begin{corollary}\label{C:closed}
Let $\sigma\in S_n$, $\tau\in S_m$, and $1\leq i\leq n$.
If $\sigma$ and $\tau$ are both closed, then there is exactly one closed
permutation $\rho$ in the interval $[B_i(\sigma,\tau),T_i(\sigma,\tau)]$.
Otherwise, this interval contains no closed permutations. 
\end{corollary}
\begin{proof} This follows from~\eqref{E:operad-M},~\eqref{E:operad-M-Q}, and~\eqref{E:D-monomial}.
\end{proof}

For instance, the permutations $\sigma=(1,2)$ and $\tau=(2,3,1)$ are closed.
We have $B_2(\sigma,\tau)=(1,3,4,2)$, $T_2(\sigma,\tau)=(3,2,4,1)$, and the
only closed permutation between these two is $(2,3,4,1)$.

Next, we discuss the analog of the filtration $\Anon^{(k)}$ of the operad $\Anon$ (Section~\ref{S:coradical}). In this case there is a stronger result: there is not only a filtration of the operad $\calQ$  but a grading.

The space $\bigoplus_{n\geq 1}\field Q_n$ carries
a Hopf algebra structure, for which the map $\sfD:\bigoplus_{n\geq 1}\field S_n\to\bigoplus_{n\geq 1}\field Q_n$ becomes a Hopf algebra surjection~\cite[Theorem 3.3]{MR}. This is the Hopf algebra of quasi-symmetric functions.
Let $k\geq 0$. The $(k+1)$-th component of the coradical filtration of this Hopf algebra has the following linear basis:
\[\{M_\setS \mid \#\setS\leq k\}\,.\]
Let $\calQ^{(k)}$ denote the corresponding sequence of spaces, i.e., the spaces with linear bases
\[\{M_\setS \mid \#\setS\leq k\,, \ \setS\subseteq[n-1]\}_{n\geq 1}\]
and $\calQ^{k}$ the sequence of spaces with linear bases
\[\{M_\setS \mid \#\setS=k\,, \ \setS\subseteq[n-1]\}_{n\geq 1}\,.\]

\begin{corollary}\label{C:operad-Q} The sequence $\{\calQ^{k}\}_{k\geq 0}$ is a grading of the non-symmetric operad $\calQ$.
\end{corollary}
\begin{proof} This follows from the fact that $\#B_i(\setS,\setT)=\#\setS+\#\setT$, which is immediate from~\eqref{E:def-B-Q}.
\end{proof}

In this sense, the non-symmetric operad $\calQ$ is graded by the cardinality of subsets.
The associated filtration is $\{\calQ^{(k)}\}_{k\geq 0}$. The map $\sfD:\Anon\to\calQ$ sends the filtration $\{\Anon^{(k)}\}_{k\geq 0}$ to the filtration $\{\calQ^{(k)}\}_{k\geq 0}$
(since a map of connected Hopf algebras preserves the coradical filtrations, or in view of~\eqref{E:D-monomial}).

Note that the space of primitive elements $\calQ^{(0)}$ is a suboperad
of $\calQ$. 
Each homogeneous component is one-dimensional,
being spanned by $M_{\emptyset_n}$, where $\emptyset_n$ denotes the empty set viewed as a subset of $[n-1]$. The operad structure is simply
\[M_{\emptyset_n}\circ_i M_{\emptyset_m}=M_{\emptyset_{n+m-1}}\,.\]

We now turn to a finer analysis of the map $\sfD:\Anon\to\calQ$. Since it is a morphism of operads (Proposition~\ref{P:operad-F-Q}), its kernel is an
operadic ideal of $\Anon$. On the other hand, recalling the coalgebra structure of
$\bigoplus_{n\geq 1}\field S_n$, it makes sense to consider the {\em Hopf kernel} of $\sfD$,
which is the space
\[\calK_* = \ker(\sfD)\cap\ker\bigl( (\id\otimes\sfD)\circ\Delta\bigr)\,. \]
(This coincides with the standard definition of Hopf kernel~\cite[Chapter 7]{Mo} in the graded connected case.) This is particularly relevant in view of the fact that there is a vector space decomposition
\[\bigoplus_{n\geq 1}\field S_n\cong \calK_*\otimes\Bigl( \bigoplus_{n\geq 1}\field Q_n\Bigr)\,.\]
(This follows from the fact that $\sfD:\bigoplus_{n\geq 1}\field S_n\to \bigoplus_{n\geq 1}\field Q_n$  is
a morphism of Hopf algebras which admits a coalgebra splitting, see~\cite[Theorem 8.1]{AS}.) 

Since $\sfD$ is degree-preserving, $\calK_*$ is a graded space. Let $\calK_n$
denote the homogeneous component of degree $n$ and let
$\calK:=\{\calK_n\}_{n\geq 1}$ denote the corresponding sequence of
spaces. We have the following surprising result.

\begin{proposition}\label{P:hopfkernel} $\calK$ is a non-symmetric suboperad of $\Anon$.
\end{proposition}
\begin{proof} Let us say that a permutation $\sigma$ is an {\em eventual identity} if there exists $k\geq 1$ such that
\[\sigma=(\ast,\ldots,\ast,1,2,\ldots k)\]
where $\ast$ stands for an arbitrary value. 
According to~\cite[Theorem 8.2]{AS}, the set 
\[\{M_\sigma \mid \sigma \text{ is {\em not} an eventual identity}\}\]
is a linear basis of $\calK_*$.  Let $\rho\in S_{n+m-1}$ and write $P_i(\rho)=(\sigma,\tau)$, where $P_i:S_{n+m-1}\to S_n\times S_m$ is the map of Section~\ref{S:weak}. In view of~\eqref{E:operad-M-P},
it suffices to show that if  $\rho$ is an eventual identity then at least one of $\sigma$ and $\tau$ are eventual identities. Assume that $\rho=(x_1,\ldots,x_{i-1},y_1,\ldots,y_m,x_{i+1},\ldots,x_n)$ is an eventual identity: 
$\rho=(\ast,\ldots,\ast,1,2,\ldots k)$. If $i+m-1\leq m+n-1-k$ then, with the notation of 
\eqref{E:def-P}, $\{1,\ldots,k\}\in C_{3,b}^\rho$. This implies that $\sigma=(\ast,\ldots,\ast,1,2,\ldots k)$. If  $i+m-1> m+n-1-k$
then $(y_1,\ldots,y_m)=(\ast,\ldots,\ast,1,2,\ldots k)$ or $(y_1,\ldots,y_m)=(l,l+1,\ldots,k)$. But $\tau=\st(y_1,\ldots,y_m)$
so in either case it is an eventual identity.
\end{proof}

\begin{remark}\label{R:compositions} We thus have a graded vector space decomposition 
\[\Anon_*\cong \calK_*\otimes\calQ_*\,,\]
in which $\calK$ is a suboperad and $\calQ$ is a quotient
of the non-symmetric operad $\Anon$. It would be interesting
to fully describe the operad structure of $\Anon$ in terms of $\calK$ and $\calQ$.
\end{remark}

We close this section by discussing descriptions of the non-symmetric operad $\calQ$ in terms of other combinatorial objects.

 A composition of $n$ is a sequence of positive integers  $\alpha=(a_1,\ldots,a_k)$ such that $a_1+\cdots+a_k=n$. There is a standard bijection between compositions of $n$ and subsets of $[n-1]$, that associates the set $\setS=\{A_1,A_2,\ldots,A_{k-1}\}$
 to $\alpha=(a_1,\ldots,a_k)$, where
 \begin{equation}\label{E:partialsum}
A_p:=a_1+\cdots+a_p\,.
\end{equation}
  In this situation, the basis element $F_\setS$ may also be denoted $F_\alpha$.
 Equation~\eqref{E:operad-F-Q} may be reformulated as follows: given compositions
 $\alpha=(a_1,\ldots,a_h)$ of $n$ and $\beta=(b_1,\ldots,b_k)$ of $m$, we have
 \[F_\alpha\circ _i F_\beta = F_{B_i(\alpha,\beta)}\]
 where
 \[B_i(\alpha,\beta):=\begin{cases}
(a_1,\ldots,a_\ell,b_1+i-1+A_\ell,b_2,\ldots,b_{k-1},b_k+A_{\ell+1}-i,a_{\ell+2},\ldots,a_h) & \text{ if }k>1\,, \\
 (a_1,\ldots,a_\ell,a_{\ell+1}+m-1,a_{\ell+2},\ldots,a_h)        & \text{ if }k=1\,,
\end{cases}\]
 where $A_p$ is as in~\eqref{E:partialsum} and $\ell$ is defined by $A_\ell<i\leq A_{\ell+1}$.
 
 A simpler description may be obtained by indexing basis elements of $\calQ$ with binary strings.  Given $\setS\subseteq[n-1]$, let
 \[\epsilon_i:=
 \begin{cases}
+ & \text{ if $i\in\setS$,} \\
 -   & \text{ if $i\notin\setS$.}
\end{cases}\]
The map that associates the sequence  $(\epsilon_1,\ldots,\epsilon_{n-1})$ to $\setS$ is a bijection between subsets of $[n-1]$ and binary strings.
With this indexing, the operad structure of $\calQ$ takes the following form:
\[F_{(\epsilon_1,\ldots,\epsilon_{n-1})}\circ_i F_{(\delta_1,\ldots,\delta_{m-1})}=
F_{(\epsilon_1,\ldots,\epsilon_{i-1},\delta_1,\ldots,\delta_{m-1},\epsilon_{i},\ldots,\epsilon_{n-1})}\,.\]
 
 \begin{remark}
The associated symmetric operad $\calS\calQ$ is the symmetric operad $MN^0$ considered by Chapoton in~\cite[Section 2]{C}.
 It can be seen that an algebra over the operad $\calQ$
  is a vector space equipped with
two associative operations $\cdot$ and $*$ satisfying
\[(a\cdot b)*c=a\cdot(b*c) \text{ \and \ }
(a* b)\cdot c=a*(b\cdot c)\,.\]
These objects
 have been  considered by Richter under the name ``Doppelalgebren''~\cite{R}
and by Pirashvili~\cite{P}.
 \end{remark}

\section{Quotient operad: planar binary trees}\label{S:trees}
 
There are combinatorial procedures for constructing a rooted planar binary tree with $n$ internal nodes from a
permutation of $[n]$, and a subset of $[n-1]$ from such a tree, which factor the
construction of the descent set of a permutation of Section~\ref{S:descent}. In this section we show that they lead to successive operad quotients of the non-symmetric associative operad. 

Let $Y_n$ be the set of rooted, planar binary trees with $n$ internal
nodes (and thus $n+1$ leaves), henceforth called simply ``trees''.
If $t\in Y_n$, we write $n:=\abs{t}$ and refer to this number as the degree of $t$.
Let $1_0$ denote the unique element of $Y_0$ (the unique tree with one leaf and no root).

Given trees $s\in Y_n$ and $t\in Y_m$, let $s\under t\in Y_{n+m}$ be the tree obtained by gluing
the root of $t$ to the rightmost branch of $s$, and 
$s/t\in Y_{n+m}$ the tree obtained by gluing the root of $s$ to the leftmost
branch of $t$. We also set $1_0\under t=t =t/1_0$.
The {\em grafting} of $s$ and $t$ is the tree
$s\vee t\in Y_{n+m+1}:=(s/\mathrm{Y})\under t= s/(\mathrm{Y}\under t)$, where $\mathrm{Y}\in Y_1$ denotes the unique tree with one internal node. Equivalently,
$s\vee t$ is obtained by drawing a new root and joining it to the roots of $s$ and $t$. For illustrations of all these operations, see~\cite[Section 1]{ASb}.

For $n\geq 1$, every tree $t\in Y_n$ has a unique decomposition $t=t_l\vee t_r$ with
$t_l\in Y_p$, $t_r\in Y_q$, $p,q\geq 0$ and $n=p+q+1$.

 We are interested in the map $\lambda:S_n\to Y_n$, as considered in~\cite[Section 2.4]{LR1}
(equivalent constructions are described in~\cite[pp.~23-24]{St86} and~\cite[Def. 9.9]{BW97}). 
The map $\lambda$ is defined recursively as follows. We set $\lambda(\id_0)=1_0$. For $n\geq 1$, let $\sigma\in S_n$
and  $j:=\sigma^{-1}(n)$.
Let $\sigma_l:=\st(\sigma_1,\ldots,\sigma_{j{-}1})$,
$\sigma_r:=\st(\sigma_{j{+}1},\ldots,\sigma_n)$, and define 
 \begin{equation}\label{E:def-lambda}
  \lambda(\sigma) := \lambda(\sigma_l)\vee\lambda(\sigma_r)\,.
 \end{equation}

Now, given $s\in Y_n$, $t\in Y_n$, and $1\leq i\leq n$, we define a tree $B_i(s,t)\in Y_{n+m-1}$ by the following recursion. Write $s=s_l\vee s_r$ and $t=t_l\vee t_r$. Let $j:=\abs{s_l}+1$. Define
\begin{equation}\label{E:def-B-Y}
B_i(s,t):=\begin{cases}
s_l\vee(s_r\circ_{i-j} t) & \text{ if $j<i$,} \\
(s_l/t_l)\vee(t_r\under s_r)   & \text{ if $j=i$,}\\
(s_l\circ_i t)\vee s_r & \text{ if $j>i$.}
\end{cases}
\end{equation}

 \begin{lemma}\label{L:lambda} For any $\sigma\in S_n$, $\tau\in S_m$, and $1\leq i\leq n$,
 \begin{equation}\label{E:lambda}
\lambda\bigl(B_i(\sigma,\tau)\bigr)=B_i\bigl(\lambda(\sigma),\lambda(\tau)\bigr)\,.
\end{equation}
\end{lemma}
\begin{proof} This may be done by induction. We omit the details.
\end{proof}

Let $\calY$ denote the sequence of spaces $\{\field Y_n\}_{n\geq 1}$. The basis element of $\field Y_n$ corresponding to a tree $t\in Y_n$ is denoted $F_t$.

 \begin{proposition}\label{P:operad-F-Y}
The sequence $\calY$ is a non-symmetric operad under the structure maps
\begin{equation}\label{E:operad-F-Y}
F_s\circ _i F_t := F_{B_i(s,t)}\,.
\end{equation}
Moreover, the map $\Lambda:\Anon\to\calY$, $\Lambda(F_\sigma):=F_{\lambda(\sigma)}$, is
a surjective morphism of non-symmetric operads.
\end{proposition}
\begin{proof} This follows from~\eqref{E:operad-F} and Lemma~\ref{L:lambda}.
\end{proof}

Let $L:Y_n\to Q_n$ be the following map. A tree $t\in Y_n$ has $n{+}1$ leaves, which we number from $1$ to $n{-}1$
left-to-right, excluding the two outermost leaves. 
Let $L(t)$ be the set of labels of those leaves that  point left. Then $\Des=\lambda\circ L$~\cite[Section 4.4]{LR1}.

\begin{corollary}\label{C:operad-A-Y-Q} The map $\sfL:\calY\to\calQ$, $\sfL(F_t):=F_{L(t)}$, is a surjective morphism of non-symmetric operads. There is a commutative diagram of surjective morphisms of non-symmetric operads
\[\xymatrix{ \Anon\ar[r]^{\Lambda}\ar@/_1pc/[rr]_{\sfD} & \calY\ar[r]^{\sfL} & \calQ }\]
\end{corollary}
\begin{proof} This follows from Propositions~\ref{P:operad-F-Q} and~\ref{P:operad-F-Y}.
\end{proof}

\medskip

We refer to~\cite[Sec. 9]{BW97} for the definition of the {\em Tamari} partial order on $Y_n$. 
The monomial basis of $\field Y_n$ is defined by  
  \begin{equation}\label{E:def-monomial-Y}
   M_s := \sum_{s\leq t} \mu(s,t) F_t\,,
 \end{equation}
where $\mu$ denotes the M\"obius function of the Tamari order. 

To describe the map $\Lambda$ on the monomial bases, we need the notion of $132$-avoiding permutation. 
A permutation $\sigma\in S_n$ meets the pattern $132$ if there is a triple of indices $1\leq i<j<k\leq n$ such that
$\sigma(i)<\sigma(k)<\sigma(j)$; i.e., if there is a $3$-letter substring $(\sigma_i,\sigma_j,\sigma_k)$ of $\sigma$ that standardizes to $(1,3,2)$. Otherwise, it is said that $\sigma$ is $132$-avoiding.

The map $\Lambda$ is given as follows~\cite[Theorem 3.1]{ASb}

\begin{equation} \label{E:Lambda-monomial}
  \Lambda( M_\sigma) = \begin{cases}
 M_{\lambda(\sigma)}& \text{ if $\sigma$ is $132$-avoiding,}\\
     0    & \text{ otherwise.}
\end{cases}
 \end{equation}

Let $\gamma:Y_n\to S_n$ be the map denoted $\textrm{Max}$ in~\cite[Definition 2.4]{LR2}. The image of $\gamma$ consists of the set of $132$-avoiding permutations~\cite[Section 1.2]{ASb}. In addition, $\lambda\circ\gamma=\id$. 

\begin{proposition}\label{P:defP-Y} Let $n,m\geq 1$ and $1\leq i\leq n$ be fixed. There is a unique map $P_i:Y_{n+m-1}\to Y_n\times Y_m$ such that
\[\xymatrix{ Y_{n+m-1}\ar[r]^{\gamma}\ar@{-->}[d]_{P_i} & S_{n+m-1}\ar[d]^{P_i}\\
Y_n\times Y_m\ar[r]_{\gamma\times\gamma} & S_n\times S_m}\]
commutes.
\end{proposition}
\begin{proof} Since $\gamma$ is injective, the claim is equivalent to the following statement: if $P_i(\rho)=(\sigma,\tau)$ and $\rho$ is $132$-avoiding, then so are $\sigma$ and $\tau$. This can be seen from~\eqref{E:def-P}.
\end{proof}

The operad structure of $\calY$ takes the following form on the monomial basis.

\begin{proposition}\label{P:operad-M-Y}
For any trees $s\in Y_n$, $t\in Y_m$, and $1\leq i\leq n$,
\begin{equation}\label{E:operad-M-P-Y}
M_s\circ _i M_t = \sum_{P_i(r)=(s,t)} M_{r}\,.
\end{equation}
\end{proposition}
\begin{proof} This follows by applying $\Lambda$ to~\eqref{E:operad-M-P} for $\sigma=\gamma(s)$ and $\tau=\gamma(t)$, in view of~\eqref{E:Lambda-monomial} and Proposition~\ref{P:defP-Y}.
\end{proof}

Next we show that each fiber of $P_i$ is an interval for the Tamari order on $Y_n$, in analogy with the situation for $S_n$~\eqref{E:operad-M} and for $Q_n$~\eqref{E:operad-M-Q}. To this end we need to recall one more map relating
permutations to trees, the map $\rho:S_n\to Y_n$ defined in~\cite[Section 2]{ASb}. The key property that relates the maps $\lambda$, $\gamma$, and $\rho$ is~\cite[Theorem 2.1]{ASb}
\begin{equation}\label{E:galois}
\sigma\leq \gamma(r)\leq \tau \iff \lambda(\sigma)\leq r  \leq \rho(\tau)\,.
\end{equation}

Let $n,m\geq 1$ and $1\leq i\leq n$.  Define a map  $T_i: Y_n\times Y_m\to Y_{n+m-1}$  by means of the second diagram below
\begin{equation}\label{E:B-T}
\xymatrix{ Y_n\times Y_m\ar[r]^{\gamma\times\gamma}\ar[d]_{B_i} & S_n\times S_m\ar[d]^{B_i}\\
Y_{n+m-1} & S_{n+m-1} \ar[l]^{\lambda} } \qquad\qquad
\xymatrix{ Y_n\times Y_m\ar[r]^{\gamma\times\gamma}\ar@{-->}[d]_{T_i} & S_n\times S_m\ar[d]^{T_i}\\
Y_{n+m-1} & S_{n+m-1} \ar[l]^{\rho} }
\end{equation}
The first diagram commutes by Lemma~\ref{L:lambda} and the fact that $\lambda\circ\gamma=\id$.

\begin{corollary}\label{C:operad-M-Y}
For any trees $s\in Y_n$, $t\in Y_m$, and $1\leq i\leq n$,
\begin{equation}\label{E:operad-M-Y}
M_s\circ _i M_t = \sum_{B_i(s,t)\leq r\leq T_i(s,t)} M_{r}\,.
\end{equation}
\end{corollary}
\begin{proof} According to the definition of $P_i$,
\[P_i(r)=(s,t) \iff \bigl(\gamma(s),\gamma(t)\bigr)= P_i\bigl(\gamma(r)\bigr)\,.\]
By Proposition~\ref{P:connection} this is  equivalent to
\[B_i\bigl(\gamma(s),\gamma(t)\bigr) \leq\gamma(r)\leq T_i\bigl(\gamma(s),\gamma(t)\bigr)\,.\]
In view of~\eqref{E:galois} this is in turn equivalent to
\[\lambda B_i\bigl(\gamma(s),\gamma(t)\bigr) \leq r\leq \rho T_i\bigl(\gamma(s),\gamma(t)\bigr)\,,\]
and by~\eqref{E:B-T}, also to 
\[B_i(s,t)\leq r\leq T_i(s,t)\,.\]
Formula~\eqref{E:operad-M-Y} is thus equivalent to~\eqref{E:operad-M-P-Y}.
\end{proof}

Unlike the case of the operad $\calQ$, the interval $[B_i(s,t),T_i(s,t)]$ does not degenerate to a point, in general. For instance,
\[ M_{\epsffile{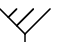}}\circ_2 M_{\epsffile{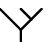}}=
 M_{\epsffile{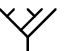}}+  M_{\epsffile{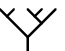}}\,.\]
Compare with~\eqref{E:ex-operad-M} and~\eqref{E:ex-operad-M-Q}.

\begin{remark}
It can be seen that an algebra over the non-symmetric operad $\calY$ is a vector space equipped with
two associative operations $\cdot$ and $*$ satisfying the identity
\[(a\cdot b)*c=a\cdot(b*c)\,.\]
This follows from results of Pirashvili~\cite{P}.
There is a different non-symmetric operad structure on the sequence $\calY$
which was introduced by Loday~\cite{Lo01} (Loday works directly with
the symmetrization of this operad). Algebras over this operad are called \emph{dendriform algebras} and they are studied in~\cite{LR1,LR2}.
The Tamari order is relevant to the structure of free dendriform algebras~\cite{LR2} .
\end{remark}

\section{The Lie operad and Dynkin's idempotent}\label{S:Lie-Dynkin}

\subsection{The associative operad as a twisted Lie algebra}\label{S:bracket}

As in~\cite[Def. 1.9]{LR2}, given permutations $\sigma\in S_n$ and $\tau\in S_m$, let
 $\sigma/\tau$ and $\sigma\under\tau$  be the following permutations in $S_{n+m}$:
 \begin{align*}
 \sigma/\tau & := ( \sigma_1,\, \sigma_2,\,\dotsc,\,\sigma_n,\ \;
     \tau_1+n,\,\tau_2+n,\,\ldots,\,\tau_m+n)\,,\\
     \sigma\under\tau &:=  ( \sigma_1+m,\, \sigma_2+m,\,\dotsc,\,\sigma_n+m,\ \,
     \tau_1,\,\tau_2,\,\ldots,\,\tau_m) \,.
 \end{align*}
The operation $F_\sigma\ast F_\tau:=F_{\sigma/\tau}$ endows $\Anon_*$ with the structure of a {\em free twisted associative algebra} in the sense of Barrat~\cite[Definitions 2 and 3]{Ba}.
Define a new operation on $\Anon_*$ by
 \begin{equation}\label{E:bracket}
\bracket{F_\sigma, F_\tau}:= F_{\sigma/\tau}-F_{\sigma\under\tau}\,.
\end{equation}
The operation $\bracket{F_\sigma, F_\tau}$ endows $\Anon_*$ with the structure of a {\em twisted Lie algebra}~\cite[Definitions 4 and 5]{Ba}.
Explicitly, this means that 
for any $x\in \Anon_n, y\in \Anon_m$, and $z\in\Anon_p$, one has
\begin{equation*}
\begin{split}
&\{y,x\}=-\{x,y\}\cdot Z_{m,n}\,, \\
&\{x,\{y,z\}\}+\{y,\{z,x\}\}\cdot Z_{n,m+p}+\{z,\{x,y\}\}\cdot Z_{n+m,p}=0\,.
\end{split}
\end{equation*}
Here, $Z_{n,m}:=1_n\under 1_m\in S_{m+m}$, and $F_{\sigma}\cdot\tau:=F_{\sigma\cdot\tau}$, where $\sigma\cdot\tau$ denotes the ordinary product of permutations in $S_n$.

It follows from~\eqref{E:operad-F} and the fact that $M_{(1,2)}=F_{(1,2)}-F_{(2,1)}$, that
 \begin{equation}\label{E:bracket-2}
\bracket{x, y}:= (M_{(1,2)}\circ_2 y)\circ_1 x
\end{equation}
for any $x,y\in\Anon$. In view of Theorem~\ref{T:operad-M}, this implies that
$\{M_\sigma,M_\tau\}$ is a linear combination of basis elements $M_\rho$
with non-negative integer coefficients. In fact, 
these structure constants are $0$ or $1$, and they admit the following description.

 Recall that an $(n,m)$-shuffle is a permutation $\zeta\in S_{n+m}$
such that
\[\zeta_1<\cdots<\zeta_n \text{ \ and \ } \zeta_{n+1}<\cdots<\zeta_{n+m}\,.\]
Let $\Sh(n,m)$ denote the set of all $(n,m)$-shuffles. This set forms an interval for the weak order of $S_{n+m}$. The smallest element is the identity $1_{n+m}:=(1,2,\ldots,n+m)$ and the biggest element is $Z_{n,m}$.



\begin{proposition}\label{P:bracket-M} For any $\sigma\in S_n$ and $\tau\in S_m$,
\begin{equation}\label{E:bracket-M}
\bracket{M_\sigma, M_\tau}=\sum_{\zeta\in\Sh(n,m) \atop \zeta\neq Z_{n,m}}
M_{\zeta\cdot(\sigma/\tau)}\,.
\end{equation}
\end{proposition}
\begin{proof} Define an operation $\ebracket'$ by means of~\eqref{E:bracket-M}.
Using~\eqref{E:fun-mon}, we calculate
\[\bracket{F_\sigma,F_\tau}'=\sum_{\sigma\leq\sigma'\,,\tau\leq\tau'}\bracket{M_\sigma,M_\tau}'= \sum_{\sigma\leq\sigma'\,,\tau\leq\tau' \atop \zeta\in\Sh(n,m)\,,\zeta\neq Z_{n,m}}
M_{\zeta\cdot(\sigma'/\tau')}\,.\]
On the other hand, according to~\cite[Prop. 2.5]{AS}, the permutations in $S_{n+m}$ that are bigger than $\sigma/\tau$ are precisely those of the form
\[\zeta\cdot(\sigma'/\tau') \text{ \ for some \ }\zeta\in\Sh(n,m),\ \sigma\leq\sigma'\in S_n\,, \text{ and }\tau\leq\tau'\in S_m\,.\]
Note that $\sigma\under\tau=Z_{n,m}\cdot(\sigma/\tau)$. Hence, again by~\cite[Prop. 2.5]{AS}, the permutations that are bigger than $\sigma\under\tau$
are those of the form
\[Z_{n,m}\cdot(\sigma'/\tau') \text{ \ for some \ } \sigma\leq\sigma'\in S_n \text{ and }\tau\leq\tau'\in S_m\,.\]
Therefore,
\begin{align*}
\bracket{F_\sigma,F_\tau}&=F_{\sigma/\tau}-F_{\sigma\under\tau}=
\sum_{\sigma\leq\sigma'\,,\tau\leq\tau' \atop \zeta\in\Sh(n,m)}M_{\zeta\cdot(\sigma'/\tau')} -
\sum_{\sigma\leq\sigma'\,,\tau\leq\tau' }M_{Z_{n,m}\cdot(\sigma'/\tau')}\\
&=\sum_{\sigma\leq\sigma'\,,\tau\leq\tau' \atop \zeta\in\Sh(n,m)\,,\zeta\neq Z_{n,m}}M_{\zeta\cdot(\sigma'/\tau')}\ =\ \bracket{F_\sigma,F_\tau}'\,.
\end{align*}
\end{proof}

\subsection{Primitive elements as a twisted Lie algebra}\label{S:primitive}
 
 Let $\Anon^{(0)}$ be the (non-symmetric) suboperad of the non-symmetric associative operad $\Anon$  discussed in Section~\ref{S:coradical}. As a vector space, $\Anon^{(0)}$ consists of
 the primitive elements of the Hopf algebra $H=\bigoplus_{n\geq 1}\field S_n$.
Therefore,  $\Anon^{(0)}_*$ is closed under the commutator bracket
\begin{equation}\label{E:comm}
[x,y]:=xy-yx\,.
\end{equation} 
Here $xy$ denotes the product in the Hopf algebra $H$ of two elements
$x,y\in H$.

This operation differs from the operation $\bracket{x,y}$ of Section~\ref{S:bracket}. For instance,
\[[F_1,F_1]=0 \text{ \ while \ }\bracket{F_1,F_1}=F_{(1,2)}-F_{(2,1)}\,.\]
In view of~\eqref{E:bracket-2} and the facts that
$M_{12}\in\Anon^{(0)}_*$ and that $\Anon^{(0)}$ is a suboperad of
$\Anon$, it follows that  $\Anon^{(0)}_*$ is  
also closed under the operation $\bracket{x,y}$. Thus $\Anon^{(0)}_*$ is a twisted Lie subalgebra of $\Anon_*$. 
More generally, if $x\in\Anon^{(k)}_*$ and
$y\in\Anon^{(h)}_*$, then $\bracket{x,y}\in\Anon^{(k+h)}_*$, by Theorem~\ref{T:coradical}.

\subsection{The Lie operad as a twisted Lie algebra}\label{S:Lie}

Let $\Lie$ be the symmetric suboperad of $\Ass$ generated by the element
\[M_{(1,2)}=F_{(1,2)}-F_{(2,1)}\,.\]
Thus $\Lie_*$ is the smallest subspace of $\Ass_*$ containing $M_{(1,2)}$, closed under the operations $\circ_i$ for every $i$, and closed under the action
of $S_n$ by right multiplication.
This is the Lie operad. Let $\calL:=\calF\Lie$ be the non-symmetric operad
obtained by forgetting the symmetric group actions (Section~\ref{S:non-symm}).


Since the suboperad $\Anon^{(0)}$ of $\Anon$ is not symmetric, we cannot immediately conclude that $\calL\subseteq \Anon^{(0)}$. Nevertheless,
Patras and Reutenauer have constructed a Hopf subalgebra of $H$ for which $\calL_*$ consists precisely of the primitive elements~\cite{PR}. It follows that
 $\calL$ is a suboperad of $\Anon^{(0)}$. The inclusions
\[\calL_n\subset\Anon^{(0)}_n\subset\Anon_n\]
are strict for $n\geq 3$.
It is known that the dimension of $\calL_n$ is $(n-1)!$~\cite[5.6.2]{Reu}. 
{}From the description of $\Anon^{(0)}$ given in Section~\ref{S:coradical} we see that
the dimension of $\Anon^{(0)}_n$ is the number of permutations with no
global descents; it lies between $(n-1)!$ and $n!$. More information on this
number is given in~\cite[Cor. 6.4]{AS}.

In view of~\eqref{E:bracket-2}, $\calL_*$ is closed under the operation $\ebracket$. Thus $\calL_*$ is a twisted Lie subalgebra of $\Anon^{(0)}_*$ and of $\Anon_*$.
Since the primitive elements of any Hopf algebra are closed under the
commutator bracket, $\calL_*$ is also closed under the commutator bracket~\eqref{E:comm}. 

\subsection{Dynkin's idempotent}\label{S:Dynkin}

The classical definition of this particular element of $\calL_n$ goes as follows~\cite[Thm. 8.16]{Reu}.
Let $V$ be a vector space. Consider the natural left action of $S_n$ on $V^{\otimes n}$ given by
\[F_\sigma\cdot v_1\ldots v_n:=v_{\sigma^{-1}(1)}\ldots v_{\sigma^{-1}(n)}\]
for any $v_1,\ldots,v_n\in V$. Let $T(V)$ be the tensor algebra of $V$,
and let $[\mu,\eta]:=\mu\eta-\eta\mu$ be the usual commutator bracket on this algebra (not to be confused with the commutator bracket on $H$ mentioned in
Section~\ref{S:primitive}).

Let $\theta_n$ be the unique element in $\field S_n$  such that for every
$V$ and $v_1,\ldots,v_n\in V$ we have
\[\theta_n\cdot(v_1\ldots v_n)=\Bigl[\cdots\bigl[ [v_1,v_2],v_3\bigr],\cdots,v_n\Bigr]\]
($n-1$ left nested commutator brackets). Dynkin's idempotent is $\frac{1}{n}\theta_n$ (it is defined when $n$ is invertible in $\field$). We may reformulate this definition as follows.

\begin{lemma}\label{L:Dynkin} For any $n\geq 1$,
\begin{equation}\label{E:Dynkin}
\theta_n=\Bigl\{\cdots\bigl\{ \{F_1,F_1\},F_1\bigr\},\cdots,F_1\Bigr\}
\end{equation}
($n-1$ left nested brackets).
\end{lemma}
\begin{proof} It suffices to show that for any $\sigma\in S_n$, $\tau\in S_m$, and
$v_1,\ldots,v_{n+m}\in V$,
\[\bracket{F_\sigma,\,F_\tau}\cdot (v_1\ldots v_{n+m})=
\bigl[F_\sigma\cdot(v_1\ldots v_n),\, F_\tau\cdot(v_{n+1}\ldots v_{n+m})\bigr]\,.\]
This follows from the facts that
\begin{align*}
F_{\sigma/\tau}\cdot (v_1\ldots v_{n+m})& =v_{\sigma^{-1}(1)}\ldots v_{\sigma^{-1}(n)}v_{n+\tau^{-1}(1)}\ldots v_{n+\tau^{-1}(m)}\\
\intertext{and}
F_{\sigma\under\tau}\cdot (v_1\ldots v_{n+m})& =v_{n+\tau^{-1}(1)}\ldots v_{n+\tau^{-1}(m)}v_{\sigma^{-1}(1)}\ldots v_{\sigma^{-1}(n)}\,.
\end{align*}
\end{proof}

We can now derive a surprisingly simple explicit expression for Dynkin's idempotent in terms of the basis $M$.
\begin{theorem}\label{T:Dynkin}
For any $n\geq 1$,
\begin{equation}\label{E:Dynkin-M}
\theta_n=\sum_{\sigma\in S_n,\,\sigma(1)=1}M_\sigma\,.
\end{equation}
\end{theorem}
\begin{proof} We argue by induction on $n$. For $n=1$, $\theta_1=F_1=M_1$.

Assume $n\geq 2$. By Lemma~\ref{L:Dynkin}, $\theta_n=\bracket{\theta_{n-1},M_1}$.
By induction hypothesis and Proposition~\ref{P:bracket-M}, we can conclude
\[\theta_n=\sum_{\tau\in S_{n-1},\,\tau(1)=1 \atop
\zeta\in\Sh(n-1,1),\, \zeta\neq Z_{n-1,1}}
M_{\zeta\cdot(\tau/1)}\,.\]
Now, $Z_{n-1,1}$ is the only $(n-1,1)$-shuffle $Z$ such that $Z(n)=1$.
All other $(n-1,1)$-shuffles $\zeta$ satisfy $\zeta(1)=1$. Thus all permutations $\sigma=\zeta\cdot(\tau/1)$ appearing in the above sum satisfy $\sigma(1)=1$.
By uniqueness of the parabolic decomposition, these $(n-1)\cdot(n-2)!$ permutations are all distinct. Hence, they are all the permutations $\sigma\in S_n$ such that $\sigma(1)=1$.
\end{proof}

Note that any permutation $\sigma\in S_n$ with $\sigma(1)=1$ has no
global descents. This confirms the fact that $\theta_n\in\Anon^{(0)}_n$.

\end{document}